\documentclass[10pt]{article}
\usepackage{amssymb}
\begin{document}
\title{The 2-color Rado Number of $x_1+x_2+\cdots +x_{m-1}=ax_m$}
\author{Dan Saracino\\Colgate University}
\date{}
\maketitle
\begin{abstract} In 1982, Beutelspacher and Brestovansky proved that for every integer $m\geq 3,$ the 2-color Rado number of the equation $$x_1+x_2+\cdots +x_{m-1}=x_m$$ is $m^2-m-1.$  In 2008, Schaal and Vestal proved that, for  every $m\geq 6,$ the 2-color Rado number of $$x_1+x_2+\cdots +x_{m-1}=2x_m$$ is $\lceil \frac{m-1}{2}\lceil\frac{m-1}{2}\rceil\rceil.$  Here we prove that, for every integer $a\geq 3$ and every $m\geq 2a^2-a+2$, the 2-color Rado number of $$x_1+x_2+\cdots +x_{m-1}=ax_m$$ is $\lceil\frac{m-1}{a}\lceil\frac{m-1}{a}\rceil\rceil.$  For the case $a=3,$ we show that our formula gives the Rado number for all $m\geq 7,$ and we determine the Rado number for all $m\geq 3$.

\end{abstract}
\vspace{.25in}

\noindent \textbf{1. Introduction}

\vspace{.25in}

A special case of the work of Richard Rado [\textbf{5}] is that for every integer $m\geq 3$ and  all positive integers $a_1,\ldots,a_m$ there exists a smallest positive integer $n$ with the following property:  for every coloring of the elements of the set $[n]=\{1,\ldots,n\}$ with two colors, there  exists a solution of the equation $$a_1x_1+a_2x_2+\cdots +a_{m-1}x_{m-1}=a_mx_m$$ using elements of $[n]$ that are all colored the same. (Such a solution is called \emph{monochromatic}.) The integer $n$ is called the \emph{2-color Rado number} of the equation.

In 1982, Beutelspacher and Brestovansky [\textbf{1}] proved that for every $m\geq 3$, the 2-color Rado number of $$x_1+x_2+\cdots +x_{m-1}=x_m$$ is $m^2-m-1.$  Since then, Rado numbers for a number of variations of this equation have also been determined.  For example, in 2008 Guo and Sun [\textbf{2}] solved the problem for the equation $$a_1x_1+a_2x_2+\cdots+a_{m-1}x_{m-1}=x_m,$$ for all positive integers $a_1,\ldots,a_{m-1}.$  They proved (confirming a conjecture of Hopkins and Schaal [\textbf{4}]) that the 2-color Rado number is $aw^2+w-a,$ where $a=\textrm{min}\{a_1,\dots, a_{m-1}\}$ and $w=a_1+\cdots +a_{m-1}.$
In the same year, Schaal and Vestal [\textbf{6}] dealt with the equation $$x_1+x_2+\cdots +x_{m-1}=2x_m.$$ They proved, in particular, that for every $m\geq 6,$ the 2-color Rado number is $\lceil\frac{m-1}{2}\lceil\frac{m-1}{2}\rceil\rceil.$  Our main purpose in the present paper is to obtain an analogue of this result for all larger values of the coefficient on $x_m.$ We prove the following result.

\vspace{.15in}

\noindent\textbf{Theorem 1}.  For every integer $a\geq 3$ and every $m\geq 2a^2-a+2,$ the 2-color Rado number of the equation$$x_1+x_2+\cdots +x_{m-1}=ax_m$$ is $\lceil\frac{m-1}{a}\lceil\frac{m-1}{a}\rceil\rceil.$

\vspace{.15in}

\noindent\textbf{Notation.}  We will denote $\lceil\frac{m-1}{a}\lceil\frac{m-1}{a}\rceil\rceil$ by $C(m,a),$ and we will denote the equation indicated in the statement of the theorem by $L(m,a).$

\vspace{.15in}

To prove Theorem 1, we show first, in Section 2,  that for all $a\geq 3$ and $m\geq 3,$ $C(m,a)$ is a lower bound for the Rado number, i.e., either $C(m,a)=1$ or there exists a 2-coloring of $ [C(m,a)-1]$ that admits no monochromatic solution of $L(m,a).$ Then, in Sections 3, 4 and 5, we show that for all $a\geq 3$ and $m\geq 2a^2-a+2,$ $C(m,a)$ is an upper bound for the Rado number, i.e., every 2-coloring of $[C(m,a)]$ admits a monochromatic solution of $L(m,a).$

In Section 6,  we prove the following.

\vspace{.15in}

\noindent \textbf{Theorem 2.}  The 2-color Rado number of $L(m,3)$ is $C(m,3)$ when $m\geq 7.$  For $m=6,5,4,3$ the Rado number is, respectively, $5,4,1,9.$

\vspace{.15in}

We confine ourselves to $m\geq 3$ because the Rado number of $L(2,a)$ fails to exist for $a\geq 3,$ as we can see by using the argument used in $[\textbf{6}]$ for $a=2.$

\vspace{.15in}     

\noindent\textbf{Notation.}  In working with a fixed 2-coloring of a set, we will use the colors red and blue, and we will denote by $R$ and $B$, respectively, the sets of elements colored red and blue. 

\vspace{.25in}

\noindent\textbf{2. Lower bounds, and some results for dealing with upper bounds}

\vspace{.25in}

\noindent\textbf{Proposition 1.}  For every $a\geq 3$ and  $m\geq 3$, the 2-color Rado number of $L(m,a)$ is at least $C(m,a).$

\vspace{.15in}

\noindent\emph{Proof.}  Since $m\geq 3$,  the Rado number exists. If $m\leq a+1$ then $C(m,a)=1$ and our claim is clear.

 Now suppose that $m\geq a+2,$ so that $C(m,a)\geq 2.$ We must show that there exists a 2-coloring of $[C(m,a)-1]$ that yields no monochromatic solution of $L(m,a).$  We use the same coloring that Schaal and Vestal used in [\textbf{5}] to establish their lower bounds, but in our less specific situation it is easier to work directly from the meaning of $C(m,a)$ than from algebraic  expressions for the ceiling function, as Schaal and Vestal did. We will use the fact that, for every real number $r$, $$r> \lceil r \rceil-1.$$

We 2-color $[C(m,a)-1]$ by coloring all the elements of $[ \lceil\frac{m-1}{a}\rceil-1]$ red and all remaining elements blue. For any red elements $x_1,x_2,\ldots,x_{m}$ we have $$\frac{x_1+\cdots +x_{m-1}}{a} \geq \frac{m-1}{a}>\left\lceil\frac{m-1}{a}\right\rceil-1\geq x_m,$$ so there are no red solutions of $L(m,a).$ For any blue elements $x_1,x_2,\ldots,x_m$ we have $$\frac{x_1+\cdots +x_{m-1}}{a}\geq \frac{m-1}{a}\left \lceil\frac{m-1}{a}\right \rceil>C(m,a)-1\geq x_m,$$ so there are no blue solutions either.   $\Box$

\vspace{.15in}

In considering upper bounds, we will often need to exhibit solutions of $L(m,a)$ in $[C(m,a)].$  To do this we will need to know that certain numbers are less than or equal to $C(m,a).$

\vspace{.15in}

 \noindent\textbf{Lemma 1.}   Suppose $a\geq 3$ and $m\geq 2a^2-a+2.$  Then the numbers $2m-2$ and $a+1$ are both less than or equal to $C(m,a).$  

\vspace{.15in}

\noindent\emph{Proof.}  By our assumption about the size of $m,$ it is clear that $2m-2\geq a+1,$ so it will suffice to prove the inequality $$2m-2\leq C(m,a).$$ 

We first consider the case $2a^2-a+2\leq m\leq 2a^2+1.$  In this case we can write $m=2a^2-a+b,$ where $2\leq b\leq a+1.$  We have $\frac{m-1}{a}=2a-1+\frac{b-1}{a},$
and therefore $\lceil \frac{m-1}{a} \rceil=2a$ and $C(m,a)=4a^2-2a+2b-2=2m-2.$

 To prove the inequality when $m\geq 2a^2+2,$ it will suffice to show that $$\frac{(m-1)^2}{a^2}\geq 2m-2.$$  This is equivalent to $$m\geq \left(2-\frac{2}{m}\right)a^2+2-\frac{1}{m},$$  and is therefore clear, since $m\geq 2a^2+2.$
 $\Box$
 
 \vspace{.15in}
 
In their treatment of upper bounds in [\textbf{5}], Schaal and Vestal proceeded by fixing the coloring of the element 1 and considering the two possibilities for the coloring of the element 2.  In dealing with $a\geq 3,$ we find it convenient to fix the coloring of the element $a-2$ and consider the two  options for the coloring of $a-1.$ 
\vspace{.15in}

\noindent\textbf{Convention.}  In  dealing with 2-colorings of $ [C(m,a)]$ in the following sections, we will assume without loss of generality that $a-2\in R.$

\vspace{.25in}

\noindent\textbf{3. Monochromatic solutions when $a-1\in B$}

\vspace{.25in}

Throughout this section we assume that $a\geq 3$ and $m\geq 2a^2-a+2.$

In this and the following sections, it will be convenient to have a compact notation for indicating solutions of $L(m,a).$

\vspace{.15in}

\noindent\textbf{Notation.}  If $n_1,\ldots,n_k$ are nonnegative integers whose sum is $m$, and $d_1,\ldots,d_k$ are elements of $[C(m,a)]$ such that we obtain a true equation from $L(m,a)$ by substituting $d_1$ for the variables $x_1,\ldots, x_{n_1}$, $d_2$ for the next $n_2$ variables, and so on, then we denote this true equation by $$[n_1\rightarrow d_1;\ n_2\rightarrow d_2;\ \cdots;\  n_k\rightarrow d_k].$$

\vspace{.15in}

For example, the true instance 
\begin{equation}
 a+a+\cdots +a=a(m-1)
\end{equation}
 of $L(m,a)$ will be denoted by $$[m-1\rightarrow a; \ 1\rightarrow m-1].$$
 
 Suppose now that we have a 2-coloring of $[C(m,a)]$ that yields no monochromatic solution of $L(m,a).$  We seek a contradiction. We will proceed by noting a number of solutions of $L(m,a)$;  all the numbers used in these solutions will be in $[C(m,a)]$ by Lemma 1.
 
 Recall that we are assuming that $a-2\in R,$ and, in this section, that $a-1\in B.$
 
 \vspace{.15in}
 
 \noindent\textbf{Lemma 2.}  We have $m-2\in R$ and $m-3\in B.$
 
 \vspace{.15in}
 
 \noindent\emph{Proof.}  Since $a-1\in B$ and there are no monochromatic solutions of $L(m,a)$ in $[C(m,a)]$, the solution
 
 \begin{equation} [m-2\rightarrow a-1; \ 2\rightarrow m-2]
 \end{equation}
  tells us that $m-2\in R.$ Likewise, since $a-2\in R,$ the solution
  \begin{equation}
  [m-3\rightarrow a-2;  \ 3\rightarrow m-3]
  \end{equation} tells us that $m-3\in B.$  $\Box$
  
  \vspace{.15in}
  
  \noindent\textbf{Lemma 3.}  We have $a\in R$ and $m-1\in B.$
  
  \vspace{.15in}
  \noindent\emph{Proof}.  Since $a-1$ and $m-3$ are in $B$, the solution
  $$ [2a\rightarrow a-1;\ m-1-2a\rightarrow a;\ 1\rightarrow m-3]
  $$ tells us that $a\in R$, and then solution (1) above tells us that $m-1\in B.$  $\Box$
  
  \vspace{.15in}
  
  Now if $m$ is even we obtain the desired contradiction by observing that the solution
  \begin{equation} \left[\frac{m-2}{2}\rightarrow a-2;\ \frac{m-2}{2}\rightarrow a;\ 2\rightarrow m-2 \right]
  \end{equation} is red.  If $m$ is odd we obtain a contradiction by considering the color of the element $a+1.$  If $a+1\in B$ then the blue solution
  \begin{equation} \left[\frac{m-1}{2}\rightarrow a-1;\ \frac{m-1}{2}\rightarrow a+1;\ 1\rightarrow m-1\right]
  \end{equation} yields a contradiction. If $a+1\in R$ then the red solution
  \begin{equation}
  \left[\frac{m-1}{2}\rightarrow a-2;\ \frac{m-5}{2}\rightarrow a;\ 1\rightarrow a+1;\ 2\rightarrow m-2\right]
  \end{equation} yields a contradiction.
  
  We have proved the following proposition.
  
  \vspace{.15in}
  
  \noindent\textbf{Proposition 2.}   If $a\geq 3, m\geq 2a^2-a+2,$ and $a-2\in R,$ then every 2-coloring of $[C(m,a)]$ with $a-1\in B$ yields a monochromatic solution of $L(m,a).$
  
  \vspace{.25in}
  
  \noindent\textbf{4.  Consequences of assuming $a-1\in R$ and no monochromatic solutions}
  
  \vspace{.25in}
  
  In this section we assume that $a\geq 3, m\geq 2a^2-a+2,$  $a-1\in R,$ and our 2-coloring of $[C(m,a)]$ yields no monochromatic solutions of $L(m,a).$  We derive some consequences that will be used in the next section.
  
  \vspace{.15in}
  
  \noindent\textbf{Lemma 4.}  Each of $m-1,m-2,m-3$ is in $B$, and $a\in R.$
  
  \vspace{.15in}
  
  \noindent\emph{Proof.}   Since $a-1$ and $a-2$ are in $R$, it follows from solutions (2) and (3) above that $m-2$ and $m-3$ are in $B$.  To conclude the proof, it will suffice, by solution (1) above, to show that $a\in R.$  But if $a\in B,$ then by solution (1) we have $m-1\in R$, and then the solution \begin{equation}
  [m-2a+1\rightarrow a-2;\ 2a-4\rightarrow a-1;\ 3\rightarrow m-1]  \end{equation} is red, a contradiction.  $\Box$
  
  \vspace{.15in}
  \noindent\textbf{Lemma 5.} We have $1\in R.$
  
  \vspace{.15in}
  \noindent\emph{Proof.}  The solution \begin{equation}
  [m-a\rightarrow 1;\ a-2\rightarrow m-1;\ 2\rightarrow m-2] \end{equation} would be blue if 1 were in $B$.  $\Box$
  
  \vspace{.15in}
  
   \noindent\textbf{Lemma 6.}  We have $2\in R$.
  
  \vspace{.15in}
 
  \noindent\emph{Proof.}  We recall that, by Lemma 1, $2m-2$ and all smaller numbers are available in $[C(m,a)]$ for use in producing solutions of $L(m,a).$
  
   From the solution $$ 
 [m-a\rightarrow 1;\ a\rightarrow m-a]  $$ we conclude that $m-a\in B$, and from the solution $$
  [m-2a+5\rightarrow a;\ a-3\rightarrow a-1;\ a-3\rightarrow 1;\ 1\rightarrow m-a+2]  $$ we conclude that $m-a+2\in B$.  If $2a\in B$, then from the solution $$[a-2\rightarrow m-a;\ m-a\rightarrow 2;\ 1\rightarrow 2a;\ 1\rightarrow m-a+2]  $$ we infer that $2\in R.$ 
  
   For the remainder of the proof we assume that $2a\in R$ and $2\in B$ and seek a contradiction.  First, by doubling all the entries in solution (1) we conclude that $2(m-1)\in B, $ and then by doubling all the entries in solution (8) we conclude that $2(m-2)\in R.$ Using the solution $$[m-2\rightarrow 2(a-1); \ 2\rightarrow 2(m-2)],$$ we see that $2(a-1)\in B,$  and then by doubling all the entries in solution (7) we see that $2(a-2)\in R.$  
   
   If $m$ is even, we get a contradiction by doubling all the entries in solution (4) to get a red solution. If $m$ is odd,  we double all the entries in solution (5) and conclude that
$2(a+1)\in R$, and then we get a contradiction by doubling all the entries in solution (6)
to get a red solution.  $\Box$
   \vspace{.15in}

 \noindent\textbf{Lemma 7.}  The numbers $m-3,m-2,m-1, \ldots, 2m-2$ are all in $B$.
 
 \vspace{.15in}
 
 \noindent\emph{Proof.}  By Lemma 4, we only need to prove this for $m,m+1,\ldots,2m-2.$ 
 
   Since $2\in R$, if $2a\in B$ we can repeat all the steps in the last two paragraphs of the proof of Lemma 6, with all the colors reversed, to obtain a contradiction. Therefore $2a\in R.$
 
 Now consider the number $m+k$, where $0\leq k\leq m-2.$ The solution  $$  [m-k-2\rightarrow a;\ k+1\rightarrow 2a;\ 1\rightarrow m+k]  $$ shows that $m+k\in B.$  $\Box$
 
 \vspace{.15in}
 
 \noindent\textbf{Lemma 8.}  The numbers $1,\ldots,2a-2$ are all in $R$.
 
 \vspace{.15in}
 
 \noindent\emph{Proof.}   We want to show that $2a-2j\in R$ for all integers $j$ such that $2\leq 2j\leq 2a-2,$ and that $2a-(2j+1)\in R$ for all $j$ such that $3\leq 2j+1\leq 2a-1.$

For $2a-2j$ we consider the solution $$ [m-(j+1)\rightarrow 2a-2j;\ j+1\rightarrow 2m-2(j+1)]  $$ and need to know that $2m-2(j+1)\in B.$  This will be true by Lemma 7 if $2m-2(j+1)\geq m-3,$ i.e., if $m\geq 2j-1.$  But this inequality holds, since $2j\leq 2a-2.$

For  $2a-(2j+1)$ we consider the solution $$ [m-(j+2)\rightarrow 2a-(2j+1);\ j\rightarrow 2m-(2j+5);\ 1\rightarrow m-2;\ 1\rightarrow2m-2(j+2)]  $$
and need to know that $2m-(2j+5)$ and $2m-2(j+2)$ are in $B$.  This will be true if
$2m-(2j+5)\geq m-3,$ i.e., if $m\geq 2j+2.$ This  inequality holds because $(2j+1)\leq 2a-1.$  $\Box$

\vspace{.15in}

There is one more result that we will need in Section 5.

\vspace{.15in}

\noindent\textbf{Lemma 9.}  If $d$ is an integer such that $a|d$ and $m-1\leq d\leq 2m-2$, then $\frac{d}{a}\in B.$

\vspace{.15in}

\noindent\emph{Proof.}  Write $d=m-1+k,$ with $0\leq k\leq m-1.$  Then the solution
$$ \left[m-1-k\rightarrow 1;\ k\rightarrow 2;\ 1\rightarrow \frac{d}{a}\right] $$ shows that $\frac{d}{a}\in B.$  $\Box$

\vspace{.25in}

\noindent\textbf{5. Monochromatic solutions when $a-1\in R$} 

\vspace{.25in}

In this section we suppose that $a\geq 3, m\geq 2a^2-a+2, a-1\in R$ and there are no monochromatic solutions of $L(m,a)$ in $[C(m,a)].$  We again seek a contradiction.

We will use the results of Section 4, and we will also need algebraic expressions for $C(m,a).$

\vspace{.15in}

\noindent\textbf{Lemma 10.}  Let $m=ua^2+va+c,$ with $u$ as large as possible and $0\leq v,c\leq a-1.$
\begin{itemize}
\item[(i)]  If $c=1$ then $C(m,a)=\frac{(m-1)^2}{a^2}.$
\item[(ii)]  If $c=0$ then $C(m,a)=\frac{m^2-m+va}{a^2}.$
\item[(iii)] If $2\leq c\leq a-1$ then $C(m,a)=\frac{m^2+(a-c-1)m+c-ac-vac+va+ta^2}{a^2},$ \\where $t=\left\lceil\frac{(c-1)(v+1)}{a}\right\rceil.$
\end{itemize}

\vspace{.15in}

\noindent\emph{Proof.}  If $c=1$ then $a|(m-1),$ so the claim is clear from the definition of $C(m,a).$  If $c=0$ then $$C(m,a)=\left\lceil\frac{m-1)}{a}\cdot \frac{m}{a}\right\rceil=\left\lceil\frac{m^2-m}{a^2}\right\rceil=\frac{m^2-m+va}{a^2},$$  since $a^2|m^2$ and $va$ is the smallest number we can add to $m^2-m$ to produce a multiple of $a^2.$

If $2\leq c\leq a-1,$ then
 $$
C(m,a)=\left\lceil\left(ua+v+\frac{c-1}{a}\right)\left (ua+v+1\right )\right\rceil$$ so $$C(m,a)=(ua+v)^2+(ua+v)+(c-1)u+t. $$
 Replacing $ua+v$ by $\left(\frac{m-c}{a}\right)$ and $u$ by $\frac{m-va-c}{a^2}$, and simplifying, we obtain the final claim of the lemma.  $\Box$

\vspace{.15in}

The three descriptions of $C(m,a)$ in Lemma 10 lead us to consider three cases.

\vspace{.15in}

\noindent  \textbf{Case 1}:  $m\equiv 1$ (mod a)

\vspace{.15in}

In this case we have $\frac{m-1}{a}\in B,$ by Lemma 9. Since $1\in R,$ we can use an idea from [\textbf{5}] and let $s$ be an integer such that $s\in R, \ s+1\in B,$ and $s+1\leq \frac{m-1}{a}.$ Then $\frac{m-1}{a}(s+1)\leq C(m,a)$ by Lemma 10,
and the solution  $$ \left[m-1\rightarrow s+1;\ 1\rightarrow \frac{m-1}{a}(s+1)\right]  $$  shows that $\frac{m-1}{a}(s+1) \in R.$  We now obtain a contradiction by noting that the solution $$ \left[\frac{m-1}{a}-a\rightarrow a-1;\ a-1\rightarrow a;\ (a-1)\frac{m-1}{a}\rightarrow s;\ 2\rightarrow \frac{m-1}{a}(s+1)\right]  $$ is red.

\vspace{.15in}

\noindent \textbf{Case 2}: $m\equiv 0$ (mod a)

\vspace{.15in}

In this case we have $\frac{m}{a}\in B$ by Lemma 9. We choose an $s$ such that $s\in R, \ s+1\in B,$ and $s+1\leq \frac{m}{a}.$  Noting that $$C(m,a)-\left(\frac{m-a}{a}\right)\frac{m}{a}=\frac{m^2-m+va}{a^2}-\frac{m^2-am}{a^2}=\frac{(a-1)m+va}{a^2},$$ we consider the element $$\alpha=\frac{m-a}{a}(s+1)+\frac{(a-1)m+va}{a^2}\leq \left(\frac{m-a}{a}\right)\frac{m}{a}+\frac{(a-1)m+va}{a^2},$$ so $\alpha\leq C(m,a).$  Noting that $\frac{m}{a}+1\in B$ by Lemma 9, we see that $\alpha\in R$ by considering the solution $$ \left[m-a\rightarrow s+1;\ v\rightarrow \frac{m}{a}+1;\ a-1-v\rightarrow \frac{m}{a};\ 1\rightarrow \alpha\right]. $$

 We now obtain a red solution of $L(m,a)$ (and therefore a contradiction) by assigning the value $\alpha$ to $x_{m-1}$ and $x_m$ and the value $s$ to $(a-1)(\frac{m-a}{a})$ other variables, and showing that we can assign values in $R$ to the remaining $\frac{m}{a}+a-3$ variables to complete the solution.  In fact we will show that we can accomplish this by using only values in the set $[2a-2]$. These values are all in $R$ by Lemma 8.

The values assigned to the remaining variables must add up to
$$\frac{a-1}{a}(m-a)+\frac{a-1}{a^2}((a-1)m+va).$$

If we can show that using only the value $2a-2$ yields a sum that is at least this large, and using only the value $1$ yields a sum that is at most this large, then there is a unique solution that uses  values  in one of the sets $\{j,\ j+1\}$, where $j\in [2a-3].$

Since $v\leq a-1,$ we can achieve our first objective by showing that 
$$(2a-2)\left(\frac{m}{a}+a-3\right)\geq \frac{a-1}{a}(m-a)+\frac{a-1}{a^2}((a-1)m+(a-1)a),$$  which simplifies to $$2a^2-8a+7-\frac{1}{a}\geq m\left(\frac{1-a}{a^2}\right),$$  and this is easily seen to be true for $a\geq 3,$ since the right-hand side is negative.

Since $v\geq 0,$ we can achieve our second objective by showing that
$$ \left(\frac{m}{a}+a-3\right)\leq \frac{a-1}{a}(m-a)+\frac{a-1}{a^2}((a-1)m).$$  But this simplifies to $2a^3-4a^2\leq m(2a^2-4a+1),$ which is true for all $a\geq 3$ and $m\ge a.$
 
\vspace{.15in}

\noindent \textbf{Case 3}: $m\equiv c$ (mod a), $2\leq c\leq a-1$

\vspace{.15in}

In this case we have $\frac{m+a-c}{a}\in B$ by Lemma 9. Choosing $s$ such that $s\in R,\ s+1\in B,$ and $s+1\leq \frac{m+a-c}{a},$ we consider the element$$
\beta=\left(\frac{m-c}{a}\right)(s+1)+\frac{(c-1)m+c-c^2-vac+va+ta^2}{a^2},$$ where $t$ is as in Lemma 10 and the second term in the sum is $$C(m,a)-\left(\frac{m-c}{a}\right)\left(\frac{m+a-c}{a}\right),$$ according to the expression for $C(m,a)$ in Lemma 10.  Then $$\beta\leq \left(\frac{m-c}{a}\right)\left(\frac{m+a-c}{a}\right)+\frac{(c-1)m+c-c^2-vac+va+ta^2}{a^2},$$ so $\beta\leq C(m,a).$  In order to work with $\beta,$ it will be helpful to have bounds on the quantity $-vac+va+ta^2.$

\vspace{.15in}

\noindent\textbf{Lemma 11.}  We have $ac-a\leq -vac+va+ta^2\leq ac-a+a^2.$

\vspace{.15in}

\noindent\emph{Proof.}  By the definition of $t$, $$\frac{(c-1)(v+1)}{a}\leq t\leq \frac{(c-1)(v+1)}{a}+1. $$  We obtain the lemma by multiplying by $a^2$ and then adding $-vac+va$.
$\Box$

\vspace{.15in}

We can now show that $\beta\in R.$  If we let  $$\gamma=\frac{(c-a)(c-1)+c-c^2-vac+va+ta^2}{a},$$
then by Lemma 11 we have $0\leq \gamma\leq a.$  By Lemma 9 we  therefore have $\frac{m+a-c}{a}+\gamma\in B,$ since $m+a-c+a^2\leq 2m-2$ because $m\geq a+a^2.$
Thus the solution $$\left[m-c\rightarrow s+1;\ c-2\rightarrow \frac{m+a-c}{a};\ 1\rightarrow \frac{m+a-c}{a}+\gamma;\ 1\rightarrow \beta\right]  $$ shows that $\beta\in R.$

To obtain our final contradiction, we construct a red solution of $L(m,a)$ by assigning the value $\beta$ to $x_{m-1}$ and $x_m$ and the value $s$ to $(a-1)(\frac{m-c}{a})$ other variables, and showing that we can assign values in $R$ to the remaining $\frac{m-c}{a}+c-2$ variables to complete the solution.  We again use values in the set $[2a-2].$

The values assigned to the remaining $\frac{m-c}{a}+c-2$ variables must add up to 
$$ \frac{a-1}{a}(m-c)+\frac{a-1}{a^2}((c-1)m+c-c^2-vac+va+ta^2),$$ which can be written as
\begin{equation}  m\left(1+\frac{c-2}{a}-\frac{c-1}{a^2}\right)+\frac{a-1}{a^2}(-ac+c-c^2-vac+va+ta^2).  \end{equation}

  If we can show that using only the value $2a-2$ (respectively, $1$) yields a sum that is at least (respectively, at most) this large, then, as before, there must be a solution that uses values in one of the sets $\{j,\ j+1\}$, where $j\in [2a-3].$ 
  
   Using the upper bound on $-vac+va+ta^2$ from Lemma 11,  we can achieve our first objective by showing that
$$(2a-2)\left(\frac{m-c}{a}+c-2\right) \geq m\left(1+\frac{c-2}{a}-\frac{c-1}{a^2}\right)+\frac{a-1}{a^2}(c-c^2-a+a^2), $$  which simplifies to $$c^2(a-1)+c(2a^3-4a^2+a+1)+(-5a^3+6a^2-a)\geq m(-a^2+1+c(a-1)).$$  If we regard $a$ as a constant and denote the quantity on the left-hand side of this inequality by $f(c)$, then the derivative $$f'(c)=2c(a-1)+(2a^3-4a^2+a+1)$$ is easily seen to be positive for $c\geq 1$ and $a\geq 3,$ so the minimum value of  $f(c)$ for $2\leq c \leq a-1$ occurs at $c=2.$ Since $$m(-a^2+1+c(a-1))\leq m(-a^2+1+(a-1)^2)=m(2-2a),$$ we only need to verify that $f(2)\geq m(2-2a),$ and this simplifies to $$2m\geq a^2+3a-2, $$ which is clearly true for $a\geq 3$ and $m\geq a^2.$

To achieve our second objective, it will suffice, by using expression (9) and the lower bound on $-vac+va+ta^2$ from Lemma 11, to show that
$$\left(\frac{m-c}{a}+c-2\right) \leq m\left(1+\frac{c-2}{a}-\frac{c-1}{a^2}\right)
+\frac{a-1}{a^2}(c-c^2-a).$$  This inequality simplifies to
$$c^2(a-1)+c(a^2-2a+1)-a^2-a\leq m(a^2-3a+1+c(a-1)).$$  Denoting the quantity on the left-hand side by $g(c)$, we have $$g'(c)=2c(a-1)+(a^2-2a+1),$$ so $g'(c)>0$ for $c\geq 1$ and $a\geq 3.$  Therefore the maximum value of $g(c)$ for $2\leq c \leq a-1$ occurs at $c=a-1.$ Since 
$$m(a^2-3a+1+c(a-1))\geq m(a^2-3a+1+2(a-1))=m(a^2-a-1),$$
we need only verify that $g(a-1)\leq m(a^2-a-1),$ i.e., that $$2a^3-7a^2+5a-2\leq m(a^2-a-1).$$  This is easily verified for $a\geq 3$ and $m\geq 2a.$

We have proved the following proposition, which completes the proof of Theorem 1.

\vspace{.15in}

\noindent\textbf{Proposition 3.}  If $a\geq 3, m\geq 2a^2-a+2,$ and $a-2\in R,$ then every 2-coloring of $[C(m,a)]$ with $a-1\in R$ yields a monochromatic solution of $L(m,a).$

\vspace{.25in}

\noindent\textbf{6. The case $a=3.$}

\vspace{.25in}

We now turn to the proof of Theorem 2.

We have determined the 2-color Rado  number of $L(m,3)$ for all $m\ge 17,$  and must consider $16\geq m\geq 3.$  We continue our convention that $a-2\in R,$ so $1\in R.$

\vspace{.15in}

\noindent\textbf{Case 1:} $16\geq m\geq 8.$

\vspace{.15in}

In this case we want to show that the Rado number is $C(m,3).$ By Proposition 1, the Rado number is at least $C(m,3).$

 We now suppose that we have a 2-coloring of $[C(m,3)]$ that yields no monochromatic solution of $L(m,3),$ and seek a contradiction. The values of $C(m,3)$ for $16\geq m\geq 8$ are, respectively, $25, 24, 22, 16,15, 14,$ $9, 8, 7. $  In each of these cases, both $a+1$ and $m-1$ are in $[C(m,a)]$, and an inspection of the arguments in Section 3 reveals that this is all we need to obtain a contradiction when $a-1\in B$, i.e., $2\in B.$  So we assume that $2\in R$, and note that then the proof of Lemma 9 is still valid, and the proof of Lemma 4 still shows that $a\in R,$ so $3\in R.$
 
 If $16\geq m\geq 14$, then $2m-8\leq C(m,3)$,  and the solution $$[m-5\rightarrow 2;\ 2\rightarrow 1;\ 3\rightarrow 2m-8]$$  shows that $2m-8\in B.$  From the solution $$[m-10\rightarrow 2;\ 8\rightarrow 1;\ 2\rightarrow m-6]$$we see that $m-6\in B,$ and then from the solution $$[m-3\rightarrow 4;\ 2\rightarrow m-6;\ 1\rightarrow 2m-8]$$we see that $4\in R.$ On the other hand, by Lemma 9, we have $5\in B$ and $6\in B.$  Now if $m=16$ then the solution $[15\rightarrow 5;\ 1\rightarrow 25]$ shows that $25\in R,$ while the solution $[ 6\rightarrow 3;\ 8\rightarrow 4;\ 2\rightarrow 25]$ shows that $25\in B.$ If $m=15$ then the solution $[12\rightarrow 5; \ 2\rightarrow 6;\  1\rightarrow 24]$ shows that $24\in R, $ while the solution $[4\rightarrow 3;\ 9\rightarrow 4;\ 2\rightarrow 24]$ shows that $24\in B.$  If $m=14,$ the solution $[12\rightarrow 5;\ 1\rightarrow 6;\  1\rightarrow 22]$ shows that $22\in R,$ while the solution $[4\rightarrow 3;\ 8\rightarrow 4;\ 2\rightarrow 22]$ shows that $22\in B.$
 
 If $13\geq m\geq 11,$ then by Lemma 9 we have $4\in B$ and $5\in B.$ If $m=13,$ then from the  solutions  $[12\rightarrow 4;\ 1\rightarrow 16]$ and  $[10\rightarrow 3;\ 1\rightarrow 2;\ 2\rightarrow 16]$ we see that $16\in R$ and $16\in B.$
 If $m=12,$ then we see from the  solutions $[10\rightarrow 4;\ 1\rightarrow 5;\ 1\rightarrow 15]$ and  $[10\rightarrow 3;\ 2\rightarrow 15]$  that $15\in R$ and $15\in B.$
 If $m=11$ then we use the solutions $[8\rightarrow 4;\ 2\rightarrow 5;\ 1\rightarrow 14]$ and $[6\rightarrow 2;\ 2\rightarrow 1;\ 3\rightarrow 14]$ to see that $14\in R$ and $14\in B.$
 
 If $10\geq m\geq 8$ then Lemma 9 shows that $3\in B,$ contradicting $3\in R.$
 
 \vspace{.15in}
 
 \noindent\textbf{Case 2:} $m=7$
 
 \vspace{.15in}
 
 In this case we again want to show that the Rado number is $C(m,3)$, which is now 4.   We know by Proposition 1 that $C(m,3)$ is a lower bound.  We no longer have have $m-1\leq C(m,a),$ however, so we cannot rely on the results of Sections 3 and 4 in showing that every 2-coloring of [4] yields a monochromatic solution of $L(m,3).$   But suppose we have a 2-coloring  that yields no such solution. Then the solution $[6\rightarrow 1;\ 1\rightarrow 2]$ shows that $2\in B,$ and doubling this solution shows that $4\in R.$  But the solution $[4\rightarrow 1;\ 3\rightarrow 4]$
 shows that $4\in B.$
 
 \vspace{.15in}
 
 \noindent\textbf{Case 3:}  $6\geq m\geq 5$
 
 \vspace{.15in} 
 
 We have $C(6,3)=4,$ and it is easy to check that the coloring $R=\{1,4\},  B=\{2,3\}$ of $[4]$ yields no monochromatic solution of $L(6,3).$  To see that the Rado number is 5, suppose we have a 2-coloring of $[5]$ that yields no monochromatic solution of $L(6,3).$ Then the solution $[4\rightarrow 1;\ 2\rightarrow 2]$ shows that $2\in B,$ and doubling this solution shows that $4\in R.$  The solution $[3\rightarrow 1;\ 3\rightarrow 3]$ shows that $3\in B.$  Then the solution $[2\rightarrow 2;\ 2\rightarrow 3;\ 2\rightarrow 5]$ shows that $5\in R,$ while the solution $[3\rightarrow 1;\ 1\rightarrow 5;\ 2\rightarrow 4]$ shows that $5\in B.$
 
 To deal with $m=5,$ note that we have $C(5,3)=3,$ and the coloring $R=\{1,3\}, B=\{2\}$ of $[3]$ yields no monochromatic solution of $L(5,3).$  To see that the Rado number is 4, suppose we have a 2-coloring of $[4]$ that yields no monochromatic solution of $L(5,3).$ The solution $[2\rightarrow 1;\ 3\rightarrow 2]$ shows that $2\in B,$ and doubling this solution shows that $4\in R.$  Then the solution $[4\rightarrow 3;\ 1\rightarrow 4]$ shows that $3\in B,$ while the solution $[3\rightarrow 2;\ 2\rightarrow 3]$ shows that $3\in R.$
 
 \vspace{.15in}
 
 \noindent\textbf{Case 4:} $4\geq m\geq 3.$
 
 \vspace{.15in}
 
 If $m=4$ then the Rado number is 1 because $[4\rightarrow 1]$ is a monochromatic solution.
 
 If $m=3$ then arguments similar to the above show that the Rado number is 9, but this result is also proved in $[\textbf{3}]$, where the 2-color Rado number of $L(3,a)$ is determined for all $a$.
 
 \vspace{.25in}

\centerline{ \textbf{References}}
\vspace{.15in}

\noindent 1.  A. Beutelspacher and W. Brestovansky, Generalized Schur numbers, \emph{Lecture Notes in Mathematics}, Springer-Verlag, Berlin, \textbf{969} (1982), 30-38

\vspace{.05in}

\noindent 2. S. Guo and Z-W. Sun, Determination of the 2-color Rado
number for $a_1x_1+ \cdots +a_mx_m=x_0$,  \emph{J. Combin. Theory Ser.A}, \textbf{115} (2008), 345-353.

\vspace{.05in}

\noindent 3.  H. Harborth and S. Maasberg,  All two-color Rado numbers for $a(x+y)=bz,$ \emph{Discrete Math.} \textbf{197-198} (1999), 397-407.

\vspace{.05in}

\noindent 4. B. Hopkins and D. Schaal, On Rado numbers for
$\displaystyle \Sigma_{i=1}^{m-1} a_ix_i=x_m$, \emph{Adv. Applied
Math.} \textbf{35} (2005), 433-441.
\vspace{.05in}

\noindent 5. R. Rado, Studien zur Kombinatorik, \emph{Mathematische
Zeitschrift} \textbf{36} (1933), 424-448.
\vspace{.05in}

\noindent 6. D. Schaal and D. Vestal, Rado numbers for $x_1+x_2+\cdots +x_{m-1}=2x_m$, \emph{Congressus Numerantium}  \textbf{191} (2008), 105-116.

\end{document}